\documentclass{elsart3-1}

\usepackage{amssymb}
\usepackage{amsmath}
\usepackage{enumerate}
\usepackage{shuffle}

\usepackage{color}

\newcommand {\bZ}{\mathbb Z}
\newcommand {\RE}{\operatorname{Re}}
\newcommand{\cH}{\mathcal{H}}
\newcommand {\bQ}{\mathbb Q}
\newcommand {\be}{\mathbf{1}}
\newcommand {\cA}{\mathcal{A}}
\def \restr#1{\mathstrut_{\textstyle |}\raise-8pt\hbox{$\scriptstyle #1$}}

\usepackage[english,francais]{babel}

\newtheorem{theorem}{Theorem}[section]
\newtheorem{lemma}[theorem]{Lemma}
\newtheorem{e-proposition}[theorem]{Proposition}

\newtheorem{e-definition}[theorem]{Definition\rm}

\newtheorem{theoreme}{Th\'eor\`eme}[section]

\newtheorem{proposition}[theoreme]{Proposition}

\renewcommand{\theequation}{\arabic{equation}}
\def\og{\leavevmode\raise.3ex\hbox{$\scriptscriptstyle\langle\!\langle$~}}
\def\fg{\leavevmode\raise.3ex\hbox{~$\!\scriptscriptstyle\,\rangle\!\rangle$}}

\begin{document}
\centerline{}
\begin{frontmatter}


\selectlanguage{english}
\title{Transfer group for renormalized\\ multiple zeta values}


\selectlanguage{english}

\author[KEF]{Kurusch Ebrahimi-Fard}
\ead{kurusch@icmat.es}
\author[DM]{Dominique Manchon}
\ead{manchon@math.univ-bpclermont.fr}
\author[JS]{Johannes Singer}
\ead{singer@math.fau.de}
\author[KEF]{Jianqiang Zhao}
\ead{zhaoj@ihes.fr}

\address[KEF]{ICMAT,C/Nicol\'as Cabrera, no.~13-15, 28049 Madrid, Spain}
\address[DM]{Univ. Blaise Pascal, C.N.R.S.-UMR 6620, 3 place Vasar\'ely, CS 60026, 63178 Aubi\`ere, France}
\address[JS]{Department Mathematik, Friedrich-Alexander-Universit\"at Erlangen-N\"urnberg, Cauerstra\ss e 11, 91058 Erlangen, Germany}

\begin{abstract}
\selectlanguage{english}
We describe in this work all solutions to the problem of renormalizing multiple zeta values at arguments of any sign in a quasi-shuffle compatible way. As a corollary we clarify the relation between different renormalizations at non-positive values appearing in the recent literature.

\vskip 0.5\baselineskip

\selectlanguage{francais}
\noindent{\bf R\'esum\'e} \vskip 0.5\baselineskip \noindent
{\bf Groupe de transfert pour les valeurs z\^eta multiples renormalis\'ees.}
Nous d\'ecrivons dans cette note l'ensemble des renormalisations possibles en des entiers de signe quelconque pour les valeurs z\^eta multiples, telles que les valeurs renormalis\'ees v\'erifient les relations de quasi-battage. En corollaire nous clarifions la relation entre les diff\'erentes renormalisations aux entiers n\'egatifs r\'ecemment apparues dans la litt\'erature.

\end{abstract}
\end{frontmatter}


\selectlanguage{francais}
\section*{Version fran\c{c}aise abr\'eg\'ee}

Les \emph{valeurs z\^eta multiples}, d\'efinies pour des entiers $k_1 \geq 2, k_2,\ldots,k_n \geq 1$ par la somme it\'er\'ee
\begin{align}\label{mzv-f}
 \zeta(k_1,\ldots,k_n):=\sum_{m_1> \cdots > m_n>0}\frac{1}{m_1^{k_1} \cdots m_n^{k_n}},
\end{align}
v\'erifient plusieurs familles de relations polynomiales, parmi lesquelles les relations de quasi-battage, qui se d\'eduisent directement de la repr\'esentation \eqref{mzv-f}, par exemple:
$$
	\zeta(a)\zeta(b)=\zeta(a,b)+\zeta(b,a)+\zeta(a+b).
$$
La d\'efinition \eqref{mzv-f} s'\'etend aux arguments complexes $(s_1,\ldots, s_n)$. La \emph{fonction z\^eta multiple} ainsi obtenue est holomorphe dans le domaine $\hbox{Re}(s_1 + \cdots + s_j) > j,\,j=1,\ldots,n$. En profondeur $n \ge 1$, la fonction z\^eta multiple s'\'etend m\'eromorphiquement  \`a $\mathbb C^n$ avec p\^oles donn\'es par~:
\begin{align*}
 s_1 & = 1, \\
 s_1 + s_2 & = 2,1,0,-2,-4,\ldots, \\
 s_1 + \cdots + s_j & \in \mathbb{Z}_{\leq j} \hbox{ for }j \geq 3.
\end{align*}
Plusieurs extensions compatibles avec les relations de quasi-battage ont \'et\'e propos\'ees, pour des arguments dans tout $\mathbb Z$  \cite{Manchon10}, ou seulement pour des arguments n\'egatifs  \cite{Guo08}, ou strictement n\'egatifs \cite{Ebrahimi15b}. Nous montrons dans cette note que toutes les solutions au probl\`eme s'obtiennent \`a partir de l'une d'entre elles via l'action libre et transitive d'un groupe de transfert, obtenu comme groupe des caract\`eres du quotient de l'alg\`ebre de Hopf des quasi-battages (sur l'alphabet $\mathbb Z$, $\mathbb Z_{\le 0}$ ou $\mathbb Z_{<0}$ selon le cas) par l'id\'eal (de Hopf) engendr\'e par les \emph{mots non-singuliers}, c'est-\`a-dire ceux associ\'es aux valeurs z\^eta multiples qui peuvent s'obtenir par prolongement analytique.\\

Nous montrons que le groupe de transfert est un groupe pro-unipotent dont l'alg\`ebre de Lie est de dimension infinie. L'espace homog\`ene principal des valeurs z\^eta multiples renormalis\'ees s'identifie donc lui aussi \`a un espace vectoriel de dimension infinie. Le point essentiel r\'eside dans la structure des mots non-singuliers, qui sont invariants par contraction de tout bloc de lettres cons\'ecutives, et qui engendrent lin\'eairement un co-id\'eal $N$ qui est de plus un co-id\'eal \`a gauche pour le coproduit de d\'econcat\'enation r\'eduit.


\selectlanguage{english}


\section{Introduction}
\label{sect:intro}

\emph{Multiple zeta values} (MZVs), which are defined for $k_1 \geq 2, k_2,\ldots,k_n\geq 1$ by the nested sum
\begin{equation}\label{mzv}
 \zeta(k_1,\ldots,k_n):=\sum_{m_1>\cdots > m_n>0}\frac{1}{m_1^{k_1}\cdots m_n^{k_n}},
\end{equation}
have a rich and deep mathematical structure \cite{Zudilin03,Hoffman05,Brown12,Zagier12}. They appeared in a 1981 preprint of J.~\'Ecalle, in the context of resurgence theory in complex analysis \cite{Ecalle81}. A systematic study however started only a decade later with the works of D.~Zagier \cite{Zagier94} and M.~E.~Hoffman \cite{Hoffman97}. M.~Kontsevich noted that \eqref{mzv} has a representation in terms of iterated integrals, which manifests itself in shuffle relations. We will focus on quasi-shuffle relations, directly derived from the nested sum representation \eqref{mzv}, for example, $\zeta(a)\zeta(b)=\zeta(a,b)+\zeta(b,a)+\zeta(a+b).$ The two representations combined yield intricate relations, which are commonly refereed to as double shuffle structures underlying MZVs \cite{Ihara06}.

Definition \eqref{mzv} extends to complex arguments $(s_1,\ldots, s_n)$. The \emph{multiple zeta function} thus obtained is holomorphic in the domain $\hbox{Re}(s_1+\cdots+s_j)>j$, $j=1,\ldots,n$. In \cite{Krattenthaler07} it was shown that the MZV function can be meromorphically continued to $\mathbb{C}^n$ with singularities
\begin{align*}
 s_1 & = 1, \\
 s_1 + s_2 & = 2,1,0,-2,-4,\ldots \\
 s_1 + \cdots + s_j & \in \bZ_{\leq j} \hbox{ ~for~}j\ge 3.
\end{align*}

This raises the natural question of how to extend MZVs to any tuple of arguments in $\mathbb Z$, such that the quasi-shuffle relations are preserved. The solution to this so-called renormalization problem for MZVs is not unique: L.~Guo and B.~Zhang gave an answer for nonpositive arguments in \cite{Guo08}. Another approach through $q$-analogues of multiple zeta values, leading to a one-parameter family of values at negative arguments, has been proposed by the first three authors in \cite{Ebrahimi15b}. An approach involving arguments in $\mathbb Z$ (including mixed signs) has been developped by S. Paycha and the second author in \cite{Manchon10}. The values obtained differ from each other. For example we have the following values for $\zeta(-1,-3)$:
\begin{align*}
 \zeta_{GZ}(-1,-3)=\frac{83}{64512}, \hspace{0.4cm}
 \zeta_{EMS}(-1,-3)= \frac{121}{94080}, \hspace{0.4cm}
 \zeta_{MP}(-1,-3)=\frac{1}{840}.
\end{align*}
Note that $\zeta_{EMS}(-1,-3)$ corresponds to the special value $t=1$ of the family
$$
	\zeta_{EMS,t}(-1,-3)= \frac{1}{8064}\frac{166t^2+166t+31}{(4t+3)(4t+1)},
$$
which is defined for $t\in \{s\in \mathbb C\colon \RE(s)>0\}$. The aforementioned solutions to the renormalization problem  of MZVs all apply the Connes--Kreimer theory \cite{Connes00} on the quasi-shuffle Hopf algebra. The different values reflect different regularization methods.

In this note, we describe how all possible quasi-shuffle compatible renormalizations for arguments in $\mathbb Z$ are related by a transfer group. More precisely, they form a set on which the transfer group acts freely and transitively. This group appears to be the character group of a commutative Hopf algebra, which is obtained as a quotient of a quasi-shuffle Hopf algebra by the (Hopf) ideal generated by \emph{non-singular words}, i.e. words associated to multiple zeta values obtainable by analytic continuation. The same approach can be used for multiple zeta values restricted to nonpositive or negative arguments. All solutions above are thus on the same orbit under the action of the corresponding restricted transfer group.


\section{Algebraic framework}

We quickly describe the quasi-shuffle Hopf algebra \cite{Hoffman00} and the set of non-singular words. Let $Y:=\{z_k\colon k\in \bZ\}$ be the set of letters. Further let $\cH:=\langle Y\rangle_\bQ$ be the free $\mathbb Q$-algebra generated by $Y$.  The empty word is denoted by $\be$. We define the \emph{quasi-shuffle product} $\ast \colon \cH\otimes \cH \to \cH$ by
\begin{enumerate}[(i)]
 \item $\be \ast w :=w\ast \be := w$,
 \item $z_m u \ast z_n v:= z_m(u \ast z_n v) + z_n(z_m u \ast v) + z_{m+n} (u \ast v)$,
\end{enumerate}
for any $w,u,v\in \cH$ and $m,n\in \bZ$. The unit map is given by $u\colon \bQ \to \cH$, $u(\lambda) = \lambda \be$. The coproduct  $\Delta\colon \cH \to \cH \otimes \cH$ is given by the deconcatenation
\begin{align}
\label{deconcat}
 \Delta(w):=\sum_{uv=w}u\otimes v
\end{align}
for any word $w\in \cH$ and the counit $\varepsilon :  \cH  \to \bQ$ is defined by $\varepsilon(\be) =1$, and $\varepsilon(w) = 0$ for any $w\neq \be$. The reduced coproduct is given by $\widetilde\Delta(w):=\Delta(w)-w\otimes\be-\be\otimes w$ for any $w\neq \be$. In Sweedler's notation we have $\widetilde\Delta(w)=\sum_{(w)}w'\otimes w''.$

Since $(\cH,\ast,\Delta)$ is a filtered and connected bialgebra it is automatically a Hopf algebra with antipode $S: \cH \to \cH$ given by $S(\be)=\be$ and
\begin{align}\label{eq:anti}
 S(w)=-w -\sum_{(w)}S(w')\ast w'' = -w -\sum_{(w)}w'\ast S(w'')
\end{align}
for any word $w\in \ker(\varepsilon)$.

A {\it{contraction of a block}} in a word $w = z_{k_1}z_{k_2}z_{k_3}\cdots z_{k_n}$ refers to combining consecutive letters in $w$ to a single letter, obtained by adding their indices together, which gives rise to a shorter word. For instance, two possible contractions in $w$ are $w_1 = z_{k_1+k_2}z_{k_3} \cdots z_{k_n}$ or $w_2 = z_{k_1}z_{k_2+k_3+k_4}z_{k_5} \cdots z_{k_n}$.

A word $w=z_{k_1}\cdots z_{k_n}$ with letters from the alphabet $Y$ will be called \emph{non-singular} if the following conditions are verified:
\begin{enumerate}[(i)]
 \item $k_1\neq 1$,
 \item $k_1+k_2\notin\{2,1,0,-2,-4,\ldots\}$,
 \item $k_1+\cdots+k_j\notin \mathbb Z_{\le j}$ for $j\ge 3$.
\end{enumerate}
Let $N \subset \cH$ be the $\bQ$-vector space spanned by the non-singular words. It is naturally graded by the length, $N=\bigoplus_{l\ge 1} N_l$. We will denote by $N_n$ the space of $\bQ$-linear combinations of non-singular words of length $n$. For any word $w=z_{k_1}\cdots z_{k_n}\in N$, the corresponding multiple zeta value
$$
	\zeta^\ast(w):=\zeta(k_1,\ldots,k_n)
$$
is either convergent or can be defined by analytic continuation \cite{Krattenthaler07}.

\begin{lemma}\label{coideal-N}
The vector space $N$ is a left coideal for the reduced deconcatenation coproduct $\widetilde{\Delta}$. Moreover, $N$ is invariant under contractions. \end{lemma}
\begin{pf}
This is straightforward and left to the reader. \qed
\end{pf}
\noindent As a corollary, it follows that $N$ is a two-sided coideal for the full deconcatenation coproduct $\Delta$ defined in \eqref{deconcat}.


\section{Renormalized multiple zeta values and the transfer group}

Let $\mathcal A$ be a commutative unital algebra, and let $\mathcal{L}(\mathcal{H},\mathcal{A})$ be the vector space of linear functions form $\mathcal{H}$ to $\mathcal{A}$. The set $G_{\mathcal A}$ of unital algebra morphisms from $\mathcal H$ to $\mathcal A$ is a group with respect to the convolution product $\star\colon \mathcal{L}(\mathcal{H},\mathcal{A})\otimes \mathcal{L}(\mathcal{H},\mathcal{A}) \to \mathcal{L}(\mathcal{H},\mathcal{A})$ defined by $ \phi \otimes \psi\mapsto  \phi \star \psi :=m_{\mathcal{A}} \circ (\phi \otimes \psi)\circ \Delta$ and with unit $e:=u\circ \varepsilon$. The inverse of $\phi \in G_{\mathcal A}$ is given by $ \phi^{-1}=\phi \circ S$.

\begin{lemma}
The set $T_{\mathcal{A}}:= \{\phi \in G_{\mathcal{A}} \colon \phi\restr{N}=0 \}$ is a subgroup of $(G_{\mathcal{A}},\star, e)$.
\end{lemma}

\begin{pf}
 Obviously, $e\in T_{\mathcal{A}}$. Let $\phi,\psi \in T_{\mathcal{A}}$. Since $G_{\mathcal{A}}$ is a group $\phi \star \psi^{-1} \in G_{\mathcal{A}}$. Further we observe for any $w\in N$ that
 \begin{align*}
  (\phi \star \psi^{-1})(w)
  & = \big(\phi \star (\psi\circ S)\big)(w) =  \phi(w) + \psi\big(S(w)\big) + \sum_{(w)}\phi(w')  \psi\big(S(w'')\big) \\
  & = \phi(w) -\psi(w) - \sum_{(w)}\psi(w') \psi\big(S(w'')\big) + \sum_{(w)}\phi(w') \psi\big(S(w'')\big) = 0
 \end{align*}
thanks to \eqref{eq:anti}, $\psi \in  G_{\mathcal{A}}$ and by means of Lemma \ref{coideal-N}.\qed
\end{pf}
The subgroup $T_{\mathcal A}$ will be called the \emph{transfer group}. Now let $\zeta:N\to \mathcal A$ be a partially defined character. We define the set of all possible renormalizations with target algebra $\cA$ by
\begin{align*}
 X_{\mathcal{A},\zeta}:= \{ \alpha \in G_{\mathcal{A}}\colon \alpha\restr{N}= \zeta \}.
\end{align*}

\begin{theorem}
The left group action $T_{\mathcal{A}} \times X_{\mathcal{A},\zeta} \to X_{\mathcal{A},\zeta}$ defined by $(\phi,\alpha)\mapsto \phi \star \alpha$ is free and transitive.
\end{theorem}

\begin{pf}
First we prove that the group action is well-defined. Let $\phi \in T_{\mathcal{A}}$ and $\alpha \in X_{\mathcal{A},\zeta}$. For any $w\in N$ we obtain
 \begin{align*}
  (\phi \star \alpha)(w) = \phi(w) + \alpha(w)  + \sum_{(w)}\phi(w') \alpha(w'') = \alpha(w) = \zeta(w)
 \end{align*}
 using Lemma \ref{coideal-N}. The identity and compatibility relations of the group action are satisfied because $X_{\mathcal{A},\zeta}\subset G_{\mathcal{A}}$.

Freeness is obvious. In order to prove transitivity, let $\alpha,\beta \in X_{\mathcal{A},\zeta}$. Then for any $w\in N$, we have
 \begin{align*}
  (\alpha \star \beta^{-1})(w)
  & = \big(\alpha \star (\beta\circ S)\big)(w) = \alpha(w) + \beta\big(S(w)\big) + \sum_{(w)}  \alpha(w')  \beta\big(S(w'')\big) \\
  & = \alpha(w) - \beta(w) - \sum_{(w)} \beta(w') \beta\big(S(w'')\big)  + \sum_{(w)} \alpha(w')  \beta\big(S(w'')\big) \\
  & = (\alpha-\beta)(w) + \sum_{(w)} (\alpha-\beta)(w')\beta\big(S(w'')\big) = 0,
 \end{align*}
since $\beta \in G_{\mathcal{A}}$ and $\alpha\restr{N}=\beta\restr{N}=\zeta$. \qed
\end{pf}
\noindent It is easily seen that the two-sided ideal $\mathcal N$ generated by $N$ is a Hopf ideal of $\cH$. The following result is straightforward.

\begin{proposition}
For any commutative unital algebra $\mathcal A$, the transfer group $T_{\mathcal A}$ is isomorphic to the group of characters of the Hopf algebra $\mathcal H/\mathcal N$.
\end{proposition}

\noindent Deconcatenation \eqref{deconcat} is conilpotent, i.e. for any word $w$ we have $\widetilde\Delta^k(w)=0$ for sufficiently large $k$, where $\widetilde \Delta^k$ is the $k$-th iterated reduced coproduct. The same property holds for the reduced coproduct of $\mathcal H/\mathcal N$. Hence, by a theorem of D. Quillen \cite{Quillen69} (see \cite[Theorem 3.9.1]{Cartier07}) the Hopf algebra $\mathcal H/\mathcal N$ is isomorphic to the free commutative algebra on a vector space $W$, the image of $\mathcal H/\mathcal N$ by the Eulerian idempotent $\pi_1=\log^\star \hbox{Id}$. As a result, the transfer group $T_{\mathcal A}$ is pro-nilpotent and can be identified with $\mathcal L(W,\mathcal A)$ by multiplicative extension. Its Lie algebra $\mathfrak t_{\mathcal A}$ is the space of $\mathcal A$-valued infinitesimal characters of $\mathcal H/\mathcal N$, which identifies with $\mathcal L(W,\mathcal A)$.

In the case $\zeta$ is given by the analytically continued multiple zeta values (with $\mathcal A=\mathbb C$), the set $X_{\mathbb C,\zeta}$ is not empty \cite{Manchon10}. There are many more elements in $X_{\mathbb C,\zeta}$, due to the fact that the transfer group $T_{\mathbb C}$ is nontrivial:

\begin{theorem}
The Lie algebra $\mathfrak t_{\mathcal A}$ of the tranfer group is infinite-dimensional.
\end{theorem}

\begin{pf}
It will be convenient to use the shuffle product instead of the quasi-shuffle. Let $\mathcal H_0$ be the Hopf algebra $\langle Y\rangle_{\mathbb Q}$ endowed with the deconcatenation and  the unital, associative, commutative shuffle product $\shuffle$ satisfying
\begin{enumerate}[(i)]
 \item $\be \shuffle w :=w \shuffle \be := w$,
 \item $z_m u \shuffle z_n v:= z_m(u \shuffle  z_n v) + z_n(z_m u \shuffle v)$.
\end{enumerate}
It is the free commutative algebra over the free Lie algebra $\hbox{Lie}(Y)$ generated by the alphabet $Y$, and is isomorphic to $\mathcal H$ by Hoffman's logarithm $\log_H$ \cite{Hoffman00}, which sends any word to a suitable linear combination of shorter words obtained by contractions. Our Hopf algebra $\mathcal H/\mathcal N$ is isomorphic to $\mathcal H_0/\mathcal N_{(0)}$, where $\mathcal N_{(0)}=\log_H(\mathcal N)$ is the ideal generated by $N_{(0)}=\log_H(N)$. Using the second assertion of Lemma \ref{coideal-N} we easily obtain $N_{(0)}=N$, and then:
$$
	W=\hbox{Lie}(Y)/\pi_1(N),
$$
where $\pi_1$ now stands for the Eulerian idempotent of $\mathcal H_0$. The natural grading on $\mathcal H_0$ by depth generates a grading on $W$. The first component $W_1$ is one-dimensional, and there are obviously infinitely many linearly independent elements in $W_2=\hbox{Lie}(Y)_2/\pi_1(N_2)$. A basis of $W_2$ can be identified with $\{[z_k,z_\ell],\, k>l \hbox{ and }z_kz_l \hbox{ is a singular word of length two}\}$. \qed
\end{pf}


\section*{Acknowledgements}
The first author is supported by a Ram\'on y Cajal research grant from the Spanish government. DM, JS and JZ would like to thank the ICMAT for its warm hospitality and gratefully acknowledge support by the Severo Ochoa Excellence Program. The second author is partially supported by Agence Nationale de la Recherche (projet CARMA).

\bibliographystyle{abbrv}
\bibliography{library}
%
%
%
%

\end{document}